\documentclass[10pt]{amsart}
\usepackage{hyperref, amssymb, amsmath, amsthm, mathrsfs}
\usepackage[numeric]{amsrefs}
\usepackage[all]{xy}

\newtheorem*{theorem*}{Theorem}

\newtheorem{theorem-def}[theorem]{Theorem-Definition}

\theoremstyle{definition}

\theoremstyle{remark}

\def    \ra     	{\rightarrow}

\def    \*      	{\times}

\def		\bZ		{\mathbb{Z}}

\def		\bZp		{\bZ_{p}}
\def		\bZpt	{\bZp^{\times}}
\def		\bG		{\mathbb{G}}
\def		\bGm	{\bG_m}

\def		\bQ		{\mathbb{Q}}
\def		\Qp		{\bQ_{p}}

\def		\bQl		{\overline{{\bQ}}_{\ell}}
\def		\bQp		{\overline{\bQ}_{p}}
\def		\bF		{\mathbb{F}}
\def		\Fp		{\bF_{p}}
\def		\Fq		{\bF_{q}}

\def		\Fp		{\mathbb{F}_{p}}

\def		\bQlt		{\bQl^\times}

\def		\fT		{\mathfrak{Tr}}

\def		\cL		{\mathcal{L}}

\def		\bW		{\mathbb{W}}

\def		\cX		{\mathcal{X}}
\def		\cY		{\mathcal{Y}}

\def		\cK		{\mathcal{K}}

\def		\Fr		{\mathrm{Fr}}

\def		\sC		{\mathsf{C}}

\newcommand         {\rar}[1]       {\stackrel{#1}{\longrightarrow}}
\newcommand         {\isom}         {\rar{\simeq}}

\newcommand         {\commentout}[1]    {}

\def		\Hom	{\operatorname{Hom}}

\newcommand{\Kum}		{\mathsf{K}}

\begin{document}

\title{Geometrization of continuous characters of $\mathbb{Z}_p^\times$}

\author{Clifton Cunningham}
\email{cunning@math.ucalgary.ca}

\author{Masoud Kamgarpour}
\email{masoudkomi@gmail.com}

\thanks{C.C. was supported by NSERC. M.K. acknowledges the hospitality and support of the University of Calgary.}


\begin{abstract} We define the $p$-adic trace of certain rank-one local systems on the multiplicative group over $p$-adic numbers, using Sekiguchi and Suwa's unification of Kummer and Artin-Schrier-Witt theories. Our main observation is that,  for every non-negative integer $n$, the $p$-adic trace defines an isomorphism of abelian groups between local systems whose order divides $(p-1)p^n$ and $\ell$-adic characters of the multiplicative group of $p$-adic integers of depth less than or equal to $n$.
\end{abstract} 

\maketitle

\commentout{

\begin{quotation}
{\smaller\smaller R\'esum\'e: Nous d\'efinissons la trace $p$-adique de certains syst\`emes
locaux de Kummer sur le groupe multiplicatif sur les nombres $p$-adiques,
\`a partir de cycles proches. Pour chaque entier positif $n$,
nous d\'emontrons que la trace $p$-adique d\'efinit un isomorphisme de
groupes ab\'eliens entre les syst\`emes locaux de Kummer dont l'ordre est
divisible par $(p-1)p^n$, et les caract\`eres $\ell$-adiques du groupe
multiplicatif des entiers $p$-adiques de profondeur plus petite ou \'egale \`a $n$.}
\end{quotation}

\vskip 20pt

}

\subsection*{Motivation}
Let $p$ and $\ell$ be distinct primes and let $q$ be a power of $p$. Let $G$ be a connected algebraic group over $\Fq$.
 To geometrize a character $\psi: G(\Fq)\ra \bQlt$ one pushes forward the Lang central extension 
\[
0\ra G(\Fq)\ra G\rar{ \operatorname{Lang}} G\ra 0,\quad \quad \operatorname{Lang}(x)=\Fr(x)-x,
\]
by $\psi^{-1}$ and obtains a local system $\cL_\psi$ on $G$. The trace of Frobenius of $\cL_\psi$ equals $\psi$; which is to say that $\cL_\psi$ and $\psi$ correspond under the functions--sheaves dictionary. Thus, we think of $\cL_\psi$ as {\it the geometrization of $\psi$}.  
Let $\sC(G)$ be the abelian group (under tensor product) consisting of $\cL_\psi$ as $\psi$ ranges over $\Hom(G(\Fq),\bQlt)$; in other words, $\sC(G)$ is the group of irreducible summands of $ \operatorname{Lang}_!\, \bQl$. Trace of Frobenius defines an isomorphism of abelian groups 
\begin{equation}\label{eq:LangIso}
t_\Fr: \sC(G)\isom \Hom(G(\Fq),\bQlt);
\end{equation}
see \cite{SGA4.5}*{Sommes Trig.} and \cite{Laumon}*{Ex. 1.1.3}. 

In this note we obtain an analogue of the above isomorphism for $\bGm$ over $p$-adic numbers. 
\begin{theorem*} 
The work of Sekiguchi and Suwa  provides an isomorphism between the abelian group of rank-one local systems on $\bG_{m,\bQp}$ whose order divides $(p-1)p^n$ and the abelian group of characters of $\bZpt$ of depth less than or equal to $n$, for every non-negative integer $n$.
\end{theorem*}
The rest of this note concerns the proof of this theorem.  



\subsection*{Unification of Kummer with Artin-Schrier-Witt}
Henceforth, we assume that $p$ is an \emph{odd} prime. Fix a non-negative integer $n$ and a primitive $p^n$-th root of unity $\zeta\in \bQp$. Set $R=\bZp[\zeta]$, $K=\Qp(\zeta)$. The main theorem of Sekiguchi and Suwa on the unification of Kummer and Artin-Schreier-Witt theories provides us with:
\begin{itemize}
\item an exact sequence 
\[
0\ra \bZ/{(p-1)\bZ}\times \bZ/{p^{n}\bZ} \ra \cY\rar{f} \cX\ra 0
\] of commutative group schemes over $R$,
\item isomorphisms $ \cY_K:=\cY\otimes_R K\isom \bG_{m,K}^{n+1}$ 
and $\cX_K \ra \bG_{m,K}^{n+1}$,
\item isomorphisms $\cY_{\Fp} \isom \bG_{m,\Fp}\times \bW_{n,\Fp}$ and $\cX_{\Fp} \isom \bG_{m,\Fp}\times \bW_{n,\Fp}$,
\end{itemize} 
such that the following diagram commutes
\[
\xymatrix@=1cm{
\ar[d]_{\theta}\bG_{m,K} & \ar[d]^{\gamma} \ar[l]_{m} \bG_{m,K}^{n+1} & \ar[r] \ar[d]^{f_K} \ar[l]_{} \cY_{K} \ar[r]  &
 \ar[d]_{f} \cY & \ar[l] \cY_{\Fp} \cong  \ar@<-40pt>[d]^{f_{\Fp}}   \bG_{m,\Fp}\times \bW_{n,\Fp} \ar@<15pt>[d]^{ \operatorname{Lang}} \\
\bG_{m,K} & \ar[l]^{\alpha} \bG_{m,K}^{n+1} &\ar[r]  \ar[l]^{} \cX_{K} \ar[r]  & \cX & \ar[l]  \cX_{\Fp} \cong \bG_{m,\Fp}\times \bW_{n,\Fp}.
}
\]
Here, $\theta(x)=x^{(p-1)p^n}$, $m$ denotes the multiplication map, $\gamma$ and $\alpha$ are defined by
\[
\gamma(x_0,...,x_{n})=(x_0^{p-1},\frac{x_1^p} {x_2}, \frac{x_2^{p}}{x_3} ,\dots,  \frac{x_{n}^p}{ x_{n-1}}),\quad 
\alpha(x_{0}, x_{1}, \ldots ,x_{n}):=\frac{(x_0 x_1 x_2 x_3 \cdots x_n)^{p^n}}{x_1 x_2^p x_3^{p^2} \cdots x_n^{p^{n-1}}}, 
\]
and $f_K$ and $f_{\Fp}$ are the restrictions of $f$ to the generic and special fibre, respectively. The main theorem of Sekiguchi and Suwa result was announced in \cite{SS-1995}. A preprint containing a proof appeared subsequently \cite{SS-1999}. According to Sekiguchi, the main tools of this preprint have been published in \cite{SS-2003}. For a general overview see \cite{Tsuchiya}.

\subsection*{The $p$-adic trace function} Let $\Kum(\bG_{m,K})$ denote the group (under tensor product) of local systems that are irreducible summands of $\theta_!\, {\bQl}$. One can easily check that all the squares in the above diagram are Cartesian; moreover, it is clear that all the vertical arrows are Galois covers of order $(p-1)p^n$. It follows that the diagram above determines a canonical isomorphism of groups 
\begin{equation}\label{eq:isoS}
S: \Kum(\bG_{m,K})\isom \sC(\bG_{m,\Fp}\times \bW_{n,\Fp}).
\end{equation}
 We define the \emph{$p$-adic trace function} by
\begin{equation} \label{eq:isoWitt}
\begin{aligned}
\fT_{n}: \Kum(\bG_{m,K}) &\longrightarrow
 \Hom(\bG_{m}(\Fp) \times \bW_{n}(\Fp),\bQlt)\\
  \cK&\mapsto t_\Fr(S(\cK)). 
  \end{aligned}
 \end{equation}
It follows at once from \eqref{eq:LangIso} and \eqref{eq:isoS} that $\fT_{n}$ is a canonical isomorphism. 

\subsection*{Relationship to continuous characters of $\bZpt$} Since $p$ is odd, the exponential map defines an isomorphism of algebraic $\Fp$-groups
\begin{equation}\label{eq:isoWitt0} 
\bG_{m,\Fp} \times \bW_{n,\Fp} \isom \bW^*_{n+1,\Fp}
\end{equation}
where $\bW^*_{n+1,\Fp}$ refers to the group scheme of units in the Witt ring scheme $\bW_{n+1,\Fp}$ (see \cite{Greenberg}) and therefore an isomorphism
\begin{equation}\label{eq:isoWitt}
\bGm(\Fp)\times \bW_n(\Fp)=\bZ/{(p-1)}\times \bZ/p^n \isom  \bZpt/(1+p^{n+1}\bZp).
\end{equation}
Accordingly, we can think of the $p$-adic trace as a character of $\bZpt/(1+p^{n+1}\bZp)$. Composing with the quotient $\bZpt \to  \bZpt/(1+p^{n+1}\bZp)$, we see that the $p$-adic trace can be interpreted as a continuous $\ell$-adic character of $\bZpt$. 
 
Conversely, for every continuous character $\chi: \bZpt\ra \bQlt$,  there is a non-negative integer $n$
such that $\chi( \bZpt/ (1+p^{n+1} \bZp) )=\{1\}$. The smallest such $n$ is known as the depth of $\chi$. We propose to think of $\cK_\chi:=\fT_{n}^{-1}(\chi)$ as {\it the geometrization of $\chi$}, when $\chi :  \bZpt\ra \bQlt$ is a continuous character of depth $n$. We do not discuss how to vary $n$ in the present text. 

We note that choosing an isomorphism of the form \eqref{eq:isoWitt} is unappetizing, to quote Deligne. 
We hope, in time, to give a construction which does not depend on this choice.

\subsection*{Relationship to character sheaves}
A character sheaf of $\bG_{m,\bQp}$ is a perverse sheaf on $\bG_{m,\bQp}$ (cohomologically) concentrated in degree $1$ where it is a rank-one local system (see \cite{Lusztig}*{\S 2}). Local systems on $\bG_{m,\bQp}$ of order dividing $(p-1)p^n$ are precisely those that have a $\Qp(\mu_{p^\infty})$-rational structure; that is, they can be defined on $\bG_{m,\Qp(\mu_{p^n})}$. In this language, the main result of this note is the following: 
{\it the $p$-adic trace of every $\Qp(\mu_{p^\infty})$-rational character sheaf on $\bG_{m,\bQp}$ is a continuous character $\bZpt\ra \bQlt$ and, moreover, every continuous $\ell$-adic character of $\bZpt$ is obtained in this manner, each one from a unique character sheaf of $\bG_{m,\bQp}$.}

\subsection*{Acknowledgement} We would like to thank T. Sekiguchi for sending us a copy of his unpublished manuscript (joint with Suwa) and for answering our questions. We thank A.-A. Aubert, J. Noel, R. Pries, T. Schedler, P. Scholze and J. Weinstein for helpful discussions and comments. Finally we would like to thank P. Deligne for carefully reading an earlier draft and providing insightful comments.


\end{document}